\documentclass[12pt,reqno,draft]{amsart}
\usepackage{amsmath,amsthm,amssymb,amscd}

\begin{document}

\newcommand{\TITLE}{Common Divisors of Elliptic Divisibility
Sequences over Function Fields}
\newcommand{\TITLERUNNING}{Common Divisors of Elliptic Divisibility Sequences}
\newcommand{\AUTHOR}{Joseph H. Silverman}
\newcommand{\DATE}{January 2004}

%%%%%%%% Set Up Theorem-Style Formats %%%%%%%%%%%%%%
\theoremstyle{plain}
\newtheorem{theorem}{Theorem}
\newtheorem{conjecture}[theorem]{Conjecture}
\newtheorem{proposition}[theorem]{Proposition}
\newtheorem{lemma}[theorem]{Lemma}
\newtheorem{corollary}[theorem]{Corollary}

\theoremstyle{definition}
% A * surpresses numbering, for example
% \newtheorem*{definition}{Definition}
\newtheorem{definition}{Definition}

\theoremstyle{remark}
\newtheorem{remark}{Remark}
\newtheorem{example}{Example}
\newtheorem{question}{Question}
\newtheorem*{acknowledgement}{Acknowledgements}

\def\BigStrut{\vphantom{$(^{(^(}_{(}$}} % Add space in tables

%%%%%%%% Set Up Environment for Notation %%%%%%%%%%%%%%
% This is currently set to allow quite wide items to be defined
\newenvironment{notation}[0]{%
  \begin{list}%
    {}%
    {\setlength{\itemindent}{0pt}
     \setlength{\labelwidth}{4\parindent}
     \setlength{\labelsep}{\parindent}
     \setlength{\leftmargin}{5\parindent}
     \setlength{\itemsep}{0pt}
     }%
   }%
  {\end{list}}

%%%%%%%% Set Up Environment for Parts in Theorems %%%%%%%%%%%%%%
\newenvironment{parts}[0]{%
  \begin{list}{}%
    {\setlength{\itemindent}{0pt}
     \setlength{\labelwidth}{1.5\parindent}
     \setlength{\labelsep}{.5\parindent}
     \setlength{\leftmargin}{2\parindent}
     \setlength{\itemsep}{0pt}
     }%
   }%
  {\end{list}}
% Use \Part{(a)}, instead of \item[(a)], to ensure upright font
\newcommand{\Part}[1]{\item[\upshape#1]}

%%%%%%%%%%%%%%%%%%
% Greek Alphabet %
%%%%%%%%%%%%%%%%%%
\renewcommand{\a}{\alpha}
\renewcommand{\b}{\beta}
\newcommand{\g}{\gamma}
\renewcommand{\d}{\delta}
\newcommand{\e}{\epsilon}
\newcommand{\f}{\phi}
\newcommand{\fhat}{{\hat\phi}}
\renewcommand{\l}{\lambda}
\renewcommand{\k}{\kappa}
\newcommand{\lhat}{\hat\lambda}
\newcommand{\m}{\mu}
\renewcommand{\o}{\omega}
\renewcommand{\r}{\rho}
\newcommand{\rbar}{{\bar\rho}}
\newcommand{\s}{\sigma}
\newcommand{\sbar}{{\bar\sigma}}
\renewcommand{\t}{\tau}
\newcommand{\z}{\zeta}

\newcommand{\D}{\Delta}
\newcommand{\F}{\Phi}
\newcommand{\G}{\Gamma}

%%%%%%%%%%%%%%%%%%%%
% Fraktur Alphabet %
%%%%%%%%%%%%%%%%%%%%
\newcommand{\ga}{{\mathfrak{a}}}
\newcommand{\gb}{{\mathfrak{b}}}
\newcommand{\gc}{{\mathfrak{c}}}
\newcommand{\gd}{{\mathfrak{d}}}
\newcommand{\gm}{{\mathfrak{m}}}
\newcommand{\gn}{{\mathfrak{n}}}
\newcommand{\gp}{{\mathfrak{p}}}
\newcommand{\gq}{{\mathfrak{q}}}
\newcommand{\gP}{{\mathfrak{P}}}
\newcommand{\gQ}{{\mathfrak{Q}}}

%%%%%%%%%%%%%%%%%%%%%%%%%
% Calligraphic Alphabet %
%%%%%%%%%%%%%%%%%%%%%%%%%
\def\Acal{{\mathcal A}}
\def\Bcal{{\mathcal B}}
\def\Ccal{{\mathcal C}}
\def\Dcal{{\mathcal D}}
\def\Ecal{{\mathcal E}}
\def\Fcal{{\mathcal F}}
\def\Gcal{{\mathcal G}}
\def\Hcal{{\mathcal H}}
\def\Ical{{\mathcal I}}
\def\Jcal{{\mathcal J}}
\def\Kcal{{\mathcal K}}
\def\Lcal{{\mathcal L}}
\def\Mcal{{\mathcal M}}
\def\Ncal{{\mathcal N}}
\def\Ocal{{\mathcal O}}
\def\Pcal{{\mathcal P}}
\def\Qcal{{\mathcal Q}}
\def\Rcal{{\mathcal R}}
\def\Scal{{\mathcal S}}
\def\Tcal{{\mathcal T}}
\def\Ucal{{\mathcal U}}
\def\Vcal{{\mathcal V}}
\def\Wcal{{\mathcal W}}
\def\Xcal{{\mathcal X}}
\def\Ycal{{\mathcal Y}}
\def\Zcal{{\mathcal Z}}

%%%%%%%%%%%%%%%%%%%%%%%%%%%%
% Blackboard Bold Alphabet %
%%%%%%%%%%%%%%%%%%%%%%%%%%%%
\renewcommand{\AA}{\mathbb{A}}
\newcommand{\BB}{\mathbb{B}}
\newcommand{\CC}{\mathbb{C}}
\newcommand{\FF}{\mathbb{F}}
\newcommand{\GG}{\mathbb{G}}
\newcommand{\NN}{\mathbb{N}}
\newcommand{\PP}{\mathbb{P}}
\newcommand{\QQ}{\mathbb{Q}}
\newcommand{\RR}{\mathbb{R}}
\newcommand{\ZZ}{\mathbb{Z}}

%%%%%%%%%%%%%%%%%%%%%%%%%%
% Boldface Math Alphabet %
%%%%%%%%%%%%%%%%%%%%%%%%%%
\def \bfa{{\mathbf a}}
\def \bfb{{\mathbf b}}
\def \bfc{{\mathbf c}}
\def \bfe{{\mathbf e}}
\def \bff{{\mathbf f}}
\def \bfF{{\mathbf F}}
\def \bfg{{\mathbf g}}
\def \bfn{{\mathbf n}}
\def \bfp{{\mathbf p}}
\def \bfr{{\mathbf r}}
\def \bfs{{\mathbf s}}
\def \bft{{\mathbf t}}
\def \bfu{{\mathbf u}}
\def \bfv{{\mathbf v}}
\def \bfw{{\mathbf w}}
\def \bfx{{\mathbf x}}
\def \bfy{{\mathbf y}}
\def \bfz{{\mathbf z}}
\def \bfX{{\mathbf X}}
\def \bfU{{\mathbf U}}
\def \bfmu{{\boldsymbol\mu}}

%%%%%%%%%%%%%%%%%%%%%%%%
% Barred Math Alphabet %
%%%%%%%%%%%%%%%%%%%%%%%%
\newcommand{\Gbar}{{\bar G}}
\newcommand{\Kbar}{{\bar K}}
\newcommand{\Obar}{{\bar O}}
\newcommand{\Pbar}{{\bar P}}
\newcommand{\Qbar}{{\bar Q}}
\newcommand{\QQbar}{{\bar{\QQ}}}

%%%%%%%%%%%%%%%%%%%%%%%%%%%%%%
% Miscellaneous New Commands %
%%%%%%%%%%%%%%%%%%%%%%%%%%%%%%

\newcommand{\Aut}{\operatorname{Aut}}
\newcommand{\Disc}{\operatorname{Disc}}
\renewcommand{\div}{\operatorname{div}}
\newcommand{\Div}{\operatorname{Div}}
\newcommand{\Etilde}{{\tilde E}}
\newcommand{\End}{\operatorname{End}}
\newcommand{\Frob}{\operatorname{Frob}}
\newcommand{\Gal}{\operatorname{Gal}}
\newcommand{\GCD}{{\operatorname{GCD}}}
\renewcommand{\gcd}{{\operatorname{gcd}}}
\newcommand{\hhat}{{\hat h}}
\newcommand{\Hom}{\operatorname{Hom}}
\newcommand{\Ideal}{\operatorname{Ideal}}
\newcommand{\Image}{\operatorname{Image}}
\newcommand{\longhookrightarrow}{\lhook\joinrel\relbar\joinrel\rightarrow}
\newcommand{\LS}[2]{\genfrac(){}{}{#1}{#2}}  % Legendre symbol
\newcommand{\MOD}[1]{~(\textup{mod}~#1)}
\newcommand{\Norm}{\operatorname{N}}
\newcommand{\NS}{\operatorname{NS}}
\newcommand{\notdivide}{\nmid}
\newcommand{\ord}{\operatorname{ord}}
\newcommand{\Pic}{\operatorname{Pic}}
\newcommand{\Proj}{\operatorname{Proj}}
\newcommand{\rank}{\operatorname{rank}}
\newcommand{\res}{\operatornamewithlimits{res}}
\newcommand{\Resultant}{\operatorname{Resultant}}
\renewcommand{\setminus}{\smallsetminus}
\newcommand{\Spec}{\operatorname{Spec}}
\newcommand{\Support}{\operatorname{Support}}
\newcommand{\tors}{{\textup{tors}}}
\newcommand{\<}{\langle}
\renewcommand{\>}{\rangle}

%%%%%%%%%%%%%%%  Topmatter %%%%%%%%%%%%%%%%%%

\title[\TITLERUNNING]{\TITLE}
\date{\DATE}
\author{\AUTHOR}
\address{Mathematics Department, Box 1917, Brown University,
Providence, RI 02912 USA}
\email{jhs@math.brown.edu}

%% \thanks{This research supported by NSA grant ***}

%%%%%%%%%%%%%%%%%%%%%%%%%%%%%%%%%%%%%%%%%%%%%%%%%%%%%%%%%%%%%%%%%%%%%%
%%%  Text (non-TeX) Abstract %%%
%% Let E/k(T) be an elliptic curve defined over a rational function field of characteristic zero. Fix a Weierstrass equation for E.  For points R in E(k(T)), write x_R=A_R/D_R^2 with relatively prime polynomials A_R(T) and D_R(T) in k[T]. The sequence {D_{nR}) for n \ge 1 is called the ``elliptic divisibility sequence of R.''
%%   Let P,Q in E(k(T)) be independent points.  We conjecture that
%%     deg (gcd(D_{nP},D_{mQ})) is bounded for m,n \ge 1,
%% and that
%%   gcd(D_{nP},D_{nQ}) = gcd(D_{P},D_{Q}) for infinitely many n \ge 1.
%% We prove these conjectures in the case that j(E) is in k. More generally, we prove analogous statements with k(T) replaced by the function field of any curve and with P and Q allowed to lie on different elliptic curves. If instead k is a finite field of characteristic p, and again assuming that j(E) is in k, we show that
%%     deg (gcd(D_{nP},D_{nQ})) >  n + O(sqrt{n})
%% for infinitely many n satisfying gcd(n,p) = 1.
%%%%%%%%%%%%%%%%%%%%%%%%%%%%%%%%%%%%%%%%%%%%%%%%%%%%%%%%%%%%%%%%%%%%%%

\begin{abstract}
Let $E/k(T)$ be an elliptic curve defined over a rational function
field of characteristic zero. Fix a Weierstrass equation for~$E$.  For
points $R\in E(k(T))$, write $x_R=A_R^{\vphantom{2}}/D_R^2$ with
relatively prime polynomials $A_R(T),D_R(T)\in k[T]$. The sequence
$\left\{D_{nR}\right\}_{n\ge1}$ is called the \emph{elliptic
divisibility sequence of~$R$}.  \par Let $P,Q\in E(k(T))$ be
independent points.  We conjecture that
\[
  \deg \bigl(\gcd(D_{nP},D_{mQ})\bigr)~\text{is bounded for $m,n\ge1$,}
\]
and that
\[
  \gcd(D_{nP},D_{nQ})
  = \gcd(D_{P},D_{Q})~\text{for infinitely many $n\ge1$.}
\]
We prove these conjectures in the case that $j(E)\in k$. More
generally, we prove analogous statements with~$k(T)$ replaced by the
function field of any curve and with~$P$ and~$Q$ allowed to lie on
different elliptic curves. If instead~$k$ is a finite field of
characteristic~$p$ and again assuming that $j(E)\in k$, we show that
$\deg \bigl(\gcd(D_{nP},D_{nQ})\bigr)$ is as large as $n+O(\sqrt{n})$
for infinitely many~$n\not\equiv0\pmod{p}$.
\end{abstract}

\maketitle

%%%%%%%%%%%%%%%%%%%%%%%%%%%%%%%%%%%%%%%%%%%%%%%%%%%%%%%%%%%%%%%%%%%%%%
\section*{Introduction}
A \textit{divisibility sequence} is a sequence~$\{d_n\}_{n\ge1}$ of
positive integers with the property that
\[
  m|n \Longrightarrow d_m|d_n.
\]
Classical examples include sequences of the form $a^n-1$ and various
other linear recurrence sequences such as the Fibonacci sequence.
See~\cite{BPvdP} for a complete classification of linear recurrence
divisibility sequences.
\par
Bugeaud, Corvaja and Zannier have shown that independent divisibility 
sequences of this type have only limited common factors. For example,
they prove that if $a,b\in\ZZ$ are  multiplicatively independent integers,
then for every $\e>0$ there is a constant $c=c(a,b,\e)$ so that
\begin{equation}
  \label{equation:BCZ}
  \log\gcd(a^n-1,b^n-1) \le \e n + c
  \qquad\text{for all $n\ge1$.}
\end{equation}
(This result is proven in~\cite{BCZ}. See also~\cite{CZ1,CZ2} for 
 more general results.)
\par
It is natural to consider the case of function fields. 
For multiplicatively independent polynomials $a,b\in k[T]$ with
coefficients in a field~$k$ of characteristic~$0$,
Ailon and Rudnick~\cite{AR} prove the strong result that there
is a constant $c=c(a,b)$ so that
\begin{equation}
  \label{equation:AR}
  \deg\gcd(a^n-1,b^n-1) \le c
  \qquad\text{for all $n\ge1$,}
\end{equation}
and that
\begin{equation}
  \label{equation:AR1}
  \gcd(a^n-1,b^n-1) =\gcd(a-1,b-1)
  \qquad\text{for infinitely many $n\ge1$,}
\end{equation}
Somewhat surprisingly, if~$a(T)$ and~$b(T)$ have coefficients in a
finite field, then neither~\eqref{equation:AR}, nor even the weaker
statement~\eqref{equation:BCZ}, is true, even if the set of
allowable~$n$'s is restricted in various reasonable ways.
(See~\cite{SilvermanGCDoverFF}.)
\par
A divisibility sequence of the form $a^n-1$ comes from a rank~1
subgroup of the multiplicative group~$\GG_m$. It is interesting to
consider divisibility sequences coming from other algebraic groups,
for example from elliptic curves. The classical definition of an
elliptic divisibility sequence~\cite{Ward1,Ward2} uses the nonlinear
relation satisfied by division polynomials, but we will use an
alternative definition\footnote{The alternative definition
of elliptic divisibility sequence that we use,
which is based directly on elliptic curves, gives a slightly different
collection of divisibility sequences than is given by the classical
non-linear recurrence formula. See~\cite[\S10.3]{EPS} and~\cite{SHIP}.}
that has the dual advantages of being more natural
and more easily generalized to other algebraic groups. (See
\cite{Durst,EGW,EPS,EW1,EW2,SHIP} for additional material on elliptic
divisibility sequences and~\cite{SilvermanGCDinFGgps} for a
discussion of general algebraic divisibility sequences and their
relation to Vojta's conjecture .)  
\par 
Let~$E/\QQ$ be an elliptic curve given by a (minimal) Weierstrass
equation
\begin{equation}
  \label{equation:WE}
  y^2+a_1xy+a_3y=x^3+a_2x^2+a_4x+a_6.
\end{equation}
Any nonzero rational point $P\in E(\QQ)$ can be written in
the form
\[
  P = (x_{P},y_{P}) 
    = \left(\frac{A_{P}}{D_{P}^2},\frac{B_{P}}{D_{P}^3}\right)
  \quad\text{with $\gcd(A_P,D_P)=\gcd(B_P,D_P)=1$.}
\]
Assume now that~$P\in E(\QQ)$ is a nontorsion point. The
\textit{elliptic divisibility sequence} associated to~$E/\QQ$ and~$P$
is the sequence of denominators of the multiples of~$P$:
\[
  \left\{ D_{nP} : n=1,2,3,\ldots\right\}
\]
The elliptic analogue of~\eqref{equation:BCZ} is part~(a) of the following 
conjecture, while part~(b) gives an elliptic analogue of 
a conjecture of Ailon and Rudnick~\cite{AR}.

\begin{conjecture}
\label{conjecture:intro:GCDonECoverQ}
With notation as above, let $P,Q\in E(\QQ)$ be independent nontorsion
points. 
\begin{parts}
\Part{(a)}
For every $\e>0$ there is a constant $c=c(E/\QQ,P,Q,\e)$ so that
\[
  \log\gcd(D_{nP},D_{nQ}) \le \e n^2 + c
  \qquad\text{for all $n\ge1$.}
\]
\Part{(b)}
There is an equality
\[
  \gcd(D_{nP},D_{nQ}) = \gcd(D_P,D_Q)
  \qquad\text{for infinitely many $n\ge1$.}
\]
\end{parts}
\end{conjecture}

\begin{remark}
Siegel's theorem~\cite[Theorem~IX.3.1]{AEC} implies that 
\[
  \log D_{nP}\gg\ll n^2,
\]
so the $n^2$ appearing in
Conjecture~\ref{conjecture:intro:GCDonECoverQ}(a) is the natural quantity
to expect. See also~\cite{SilvermanGCDinFGgps} for a proof that
Vojta's conjecture~\cite{Vojta} implies
Conjecture~\ref{conjecture:intro:GCDonECoverQ}(a).
\end{remark}

Following the lead of Ailon and Rudnick, we replace~$\QQ$ with the
rational function field~$k(T)$ and replace~$\ZZ$ with the ring of
polynomials~$k[T]$. Then we can look at an elliptic curve $E/k(T)$
given by a (minimal) Weierstrass equation~\eqref{equation:WE} and we
can write the $x$-coordinate of a point $P\in E(k(T))$ in the form
\[
  x_P=\frac{A_P}{D_P^2}
  \quad\text{with $A_P(T),D_P(T)\in k[T]$ satisfying
  $\gcd(A_P,D_P)=1$.}
\]
This leads to a conjectural elliptic analogue of~\eqref{equation:AR}
and~\eqref{equation:AR1}.

\begin{conjecture}
\label{conjecture:intro:GCDonECoverFF}
Let~$k$ be an algebraically closed field of characteristic~$0$, let
$E/k(T)$ be an elliptic curve, and let $P,Q\in E(k(T))$ be independent
points. Then there is a constant $c=c(E,P,Q)$ so that
\[
  \deg\gcd(D_{nP},D_{nQ}) \le c
  \qquad\text{for all $n\ge1$.}
\]
Further, there is an equality
\[
  \gcd(D_{nP},D_{nQ}) = \gcd(D_{P},D_{Q}) 
  \qquad\text{for infinitely many $n\ge1$.}
\]
\end{conjecture}

Our principal result in this paper is a proof of
Conjecture~\ref{conjecture:intro:GCDonECoverFF} in the case that~$E$
has constant $j$-invariant. We note that even in this special case,
the proof requires nontrivial tools such as Raynaud's
theorem~\cite{Raynaud1,Raynaud2} bounding torsion points on
subvarieties of abelian varieties.

\begin{theorem}
\label{theorem:intro:maintheorem}
Conjecture~\ref{conjecture:intro:GCDonECoverFF} is true for elliptic
curves with constant $j$-invariant, i.e., with $j(E)\in
k$.\end{theorem}

We actually prove something more general than
Theorem~\ref{theorem:intro:maintheorem}. First, we  replace~$k(T)$
by the function field of an arbitrary algebraic curve. Second, we 
allow different integer multipliers for~$P$ and~$Q$.  Third, we 
allow the points~$P$ and~$Q$ to lie on different elliptic curves. For
the complete statement, see Conjecture~\ref{conjecture:ECoverFF} and
Theorem~\ref{theorem:ECoverFFconstantj}.  This added generality does
not  significantly lengthen the proof and makes parts of the argument
more transparent.

We also consider the case that~$E$ is an elliptic curve over a field
$\FF_q(T)$ of characteristic~$p$. In this case, nothing like
Conjecture~\ref{conjecture:intro:GCDonECoverFF} is true, even with the
natural restriction that~$n$ be prime to~$p$. We prove
(Theorem~\ref{theorem:gcdonECoverFFconstantj}) that if~$E/\FF_q(T)$ has
constant~$j$ invariant and $P,Q\in E(\FF_q(T))$, then
\[
  \deg\gcd(D_{nP},D_{nQ}) \ge n + O(\sqrt{n})
\]
for infinitely many~$n$ satisfying $p\notdivide n$. We conjecture that
the same is true for all~$E/\FF_q(T)$.

\begin{acknowledgement}
The author would like to thank Gary Walsh for rekindling his interest
in the arithmetic properties of divisibility sequences and for bringing
to his attention the articles~\cite{AR} and~\cite{BCZ}.
\end{acknowledgement}

%%%%%%%%%%%%%%%%%%%%%%%%%%%%%%%%%%%%%%%%%%%%%%%%%%%%%%%%%%%%%%%%%%%%%%

\section{Preliminaries}
\label{section:preliminaries}
In this section we set some notation, recall a deep theorem,
and prove two basic estimates that will be required for our main results.
We begin with notation.

\begin{notation}
\item[$k$] 
an algebraically closed field of characteristic zero.
\item[$C/k$]
a smooth projective curve.
\item[$K$]
the function field of $C$.
\item[$E/K$]
an elliptic curve.
\item[$\Ecal/C$]
a minimal smooth projective elliptic surface $\Ecal\to C$ with generic
fiber~$E$.
\item[$\s_P$]
the section $\s_P:C\to\Ecal$ corresponding to a point $P\in E(K)$.
%% \item[$\Pbar$]
%% the divisor $\Pbar=\s_P(C)\in\Div(\Ecal)$.
\item[$\Obar$]
the ``zero divisor'' $\Obar=\s_O(C)\in\Div(\Ecal)$ corresponding to
the point $O\in E(K)$.
\end{notation}

We recall that~$E/K$ is said to \textit{split over~$K/k$} if it
is~$K$-isomorphic to an elliptic curve defined over~$k$, and
that~$E/K$ is said to be \textit{constant over~$k$} if $j(E)\in k$.
Clearly split curves are constant, while a constant curve can always
be split over a finite extension of~$K$.
\par
We observe that if a Weierstrass equation is chosen for~$E/K$ and if
$P=(x_P,y_P)\in E(K)$, then the pullback divisor $\s_P^*(\Obar)$ is,
roughly, one half the polar divisor of~$x_P$. The following elementary
result shows the stability of $\s_{mP}^*(\Obar)$ at at a fixed point
of~$C$ for multiples~$mP$ of~$P$. This is well known, but for
completeness and since it is false when the residue characteristic is
positive, we include a proof.

\begin{lemma}
\label{lemma:stabilityofvanishing}
With notation as above, let $\g\in C$ and let $P\in E(K)$ be a
nontorsion point.
\begin{parts}
\Part{(a)}
If $\ord_\g \s_{P}^*(\Obar)\ge 1$, then
\[
  \ord_\g\s_{mP}^*(\Obar) = \ord_\g\s_{P}^*(\Obar)
  \qquad\text{for all $m\ne 0$.}
\]
\Part{(b)}
There is an integer $m'=m'(E/K,P,\g)$ so that
\[
  \ord_\g\s_{mP}^*(\Obar) \in \{0,m'\}
  \qquad\text{for all $m\ne 0$.}
\]
In particular, $\ord_\g \s_{mP}^*(\Obar)$ is bounded independently of~$m$.
\end{parts}
\end{lemma}
\begin{proof} (a)
Let $[m]:\Ecal\to\Ecal$ be the multiplication-by-$m$ map. Notice that
\[
  [m]^*\Obar = \Obar + D_m,
\]
where the divisor $D_m\in\Div(\Ecal)$ is the divisor of nonzero
$m$-torsion points. For example, if the fiber~$\Ecal_\g$ is nonsingular,
then the intersection of~$D_m$ with~$\Ecal_\g$ consists of the nonzero
$m$-torsion points of the elliptic curve~$\Ecal_\g$. It is thus clear
that at least on the nonsingular fibers, the divisors~$\Obar$ and~$D_m$ do
not intersect. (This is where we are using the characteristic zero
assumption. More generally, it is enough to assume that~$m$ is
relatively prime to the residue characteristic.) However, even if
the fiber~$\Ecal_\g$ is singular, the map~$[m]:\Ecal\to\Ecal$ is \'etale
in a neighborhood at the zero point~$O_\g$ of~$\Ecal_\g$, 
so $\Obar\cap D_m=\emptyset$.
\par 
We have
\[
  \s_{mP}^*(\Obar) 
  = \s_{P}^*\bigl([m]^*(\Obar)\bigr) 
  = \s_{P}^*(\Obar) + \s_{P}^*(D_m).
\]
The assumption that $\ord_\g\s_{P}^*(\Obar)\ge1$ is equivalent to the
statement that $\s_{P}(\g)=O_\g\in\Ecal_\g$. 
Since $\Obar\cap D_m=\emptyset$ from above, it follows that
the support of $\s_{P}^*(D_m)$ does
not contain~$\g$, so
\[
  \ord_\g \s_{mP}^*(\Obar) 
  = \ord_\g\s_{P}^*(\Obar) + \ord_\g\s_{P}^*(D_m)
  = \ord_\g\s_{P}^*(\Obar).
\]
This completes the proof of~(a).  \par In order to prove~(b), we may
suppose without loss of generality that there exists some $m\ne0$ such
that $\ord_\g \s_{mP}^*(\Obar)\ge1$, since otherwise we may take
$m'=0$. Suppose that \text{$\ord_\g \s_{m_1P}^*(\Obar)\ge1$} and
\text{$\ord_\g \s_{m_2P}^*(\Obar)\ge1$}. Then applying~(a), first to~$m_1P$
with $m=m_2$ and second to~$m_2P$ with $m=m_1$, we find that
\[
  \ord_\g \s_{m_1P}^*(\Obar) =\ord_\g \s_{m_2m_1P}^*(\Obar) 
  =\ord_\g \s_{m_2P}^*(\Obar) .
\]
\end{proof}

\begin{remark}
Lemma~\ref{lemma:stabilityofvanishing} readily generalizes to algebraic
groups. For the group~$\GG_m/k(\PP^1)$, the proof is especially
transparent and helps to illustrate the general case, so we recall it
here. Let $R(T)\in k(T)$ be a rational function, say $R(T)=A(T)/B(T)$,
and suppose that $\ord_\g(R(T)-1)=e\ge 1$. This implies that $B(\g)\ne0$
and that
\[
  A(T)-B(T) = (T-\g)^e C(T)
  \qquad\text{for some $C(T)\in k[T]$.}
\]
Then
\begin{align*}
  A(T)^{m}-B(T)^{m}
  &=\prod_{\z\in\bfmu_m} \bigl(A(T)-\z B(T)\bigr)\\
  &=\prod_{\z\in\bfmu_m} \bigl((T-\g)^e C(T)+ (1-\z)B(T)\bigr).
\end{align*}
Hence
\[
  \left.\frac{A(T)^{m}-B(T)^{m}}{A(T)-B(T)}\right|_{T=\g}
  = \prod_{\substack{\z\in\bfmu_m\\\z\ne1\\}} (1-\z)B(\g)
  \ne 0,
\]
which proves that
$\ord_\g\bigl(A(T)^{m}-B(T)^{m}\bigr)=\ord_\g\bigl(A(T)-B(T)\bigr)$.
\end{remark}

We will also need the following elementary result, which says that a
$\Kbar$-isogeny mapping even one $K$-rational nontorsion point to a
$K$-rational point is necessarily itself defined over~$K$.

\begin{lemma}
\label{lemma:isogenyofrationalpt}
Let~$K$ be a field of characteristic~$0$, let~$E_1/K$ and~$E_2/K$ be
elliptic curves, and let $G:E_2\to E_1$ be an isogeny defined
over~$\Kbar$.  Suppose that there is a $K$-rational point $P\in
E_2(K)$ so that the image~$G(P)$ is also $K$-rational, i.e., $G(P)\in
E_1(K)$. Then either~$P$ has finite order or else~$G$ is defined over~$K$.
\end{lemma}
\begin{proof}
For each $s\in\Gal(\Kbar/K)$, define an isogeny~$g_s$ by
\[
  g_s:E_2\longrightarrow E_1,\qquad
  g_s(Q) = G^s(Q) - G(Q).
\]
The assumption on the point~$P$ implies that
\[
  G(P) = \bigl(G(P)\bigr)^s = G^s(P^s) = G^s(P),
\]
so we see that $P\in\ker(g_s)$ for all $s\in\Gal(\Kbar/K)$.
Let $d_s=\deg(g_s)$. Applying the dual isogeny, it follows that
$P\in E_2[d_s]$ for all $s\in\Gal(\Kbar/K)$.
Hence either~$P$ is a torsion point, or else $d_s=0$ 
for all $s\in\Gal(\Kbar/K)$.
But
\[
  d_s=0 \iff g_s=[0] \iff G^s=G,
\]
so 
\[
 \text{$d_s=0$ for all $s\in\Gal(\Kbar/K)$}
 \iff
 \text{$G$ is defined over~$K$.}
\]
This completes the proof that either~$P$ is a torsion point or else~$G$
is defined over~$K$.
\end{proof}

We conclude this section by recalling a famous result of
Raynaud. We will apply Raynaud's theorem to a curve embedding in an
abelian surface.

\begin{theorem}[Raynaud's Theorem]
\label{theorem:raynaud}
Let $k$ be a field of characteristic zero, let~$A/k$ be an abelian
variety, and let $V\subset A$ be a subvariety. Then the Zariski
closure of $V\cap A_\tors$ is equal to a finite union of translates of
abelian subvarieties of~$A$ by torsion points. 
\end{theorem}
\begin{proof}
See~\cite{Raynaud2} for the general case.  For the case that~$V$ is a
curve, which is the case that we will need, see~\cite{Raynaud1}.
\end{proof}

%%%%%%%%%%%%%%%%%%%%%%%%%%%%%%%%%%%%%%%%%%%%%%%%%%%%%%%%%%%%%%%%%%%%%%
\section{Common Divisors on Elliptic Curves over 
Characteristic~$0$ Function Fields}

We continue with the notation set in Section~\ref{section:preliminaries}.
For any two \underbar{effective} divisors $D_1,D_2\in\Div(C)$, we define
the \textit{greatest common divisor} in the usual way as
\[
  \GCD(D_1,D_2) 
  = \sum_{\g\in C} \min\bigl\{\ord_\g(D_1),\ord_\g(D_2)\bigr\}\cdot(\g)
  \in \Div(C).
\]
(Here~$\ord_\g(D)$ is the coefficient of~$\g$ in the divisor~$D$.)

\begin{definition}
Let~$E/K$ be an elliptic curve over a function field as above, and let
$P,Q\in E(K)$ be points, not both zero.  Then the
\textit{{\upshape(}elliptic{\upshape)} greatest common divisor} of~$P$
and~$Q$ is
\[
  \GCD(P,Q) = \GCD\bigl(\s_P^*(\Obar),\s_Q^*(\Obar)\bigr).
\]
\end{definition}

\begin{remark}
As noted earlier, the divisor $\s_P^*(\Obar)$ is, roughly, one half
the polar divisor of~$x_P$.  Thus the elliptic~$\GCD$ is a natural
generalization of the definition given in the introduction. 
\end{remark}

More generally, we can work with points on different curves.  

\begin{definition}
Let~$E_1/K$ and~$E_2/K$ be elliptic curves as above, and let $P_1\in
E_1(K)$ and $P_2\in E_2(K)$ be points, not both zero.
The \textit{{\upshape(}elliptic{\upshape)} greatest common divisor}
of~$P_1$ and~$P_2$ is the divisor
\[
  \GCD(P_1,P_2) = \GCD\bigl(\s_{P_1}^*(\Obar_{\Ecal_1}),
    \s_{P_2}^*(\Obar_{\Ecal_2})\bigr)\in\Div(C).
\]
\end{definition}

In order to prove boundedness of~$\GCD(P_1,P_2)$, we need~$P_1$ and~$P_2$
to be independent in some appropriate sense, which prompts the following
definition.

\begin{definition}
\label{definition:independence}
We say that~$P_1$ and~$P_2$ are \textit{dependent} if there are
isogenies $F:E_1\to E_1$ and $G:E_2\to E_1$, not both zero, so that
$F(P_1)=G(P_2)$; otherwise we say that~$P_1$ and~$P_2$ are
\textit{independent}. If~$E_1$ and~$E_2$ are defined over a field~$K$,
we say that~$P_1$ and~$P_2$ are \textit{$K$-dependent} if the
isogenies~$F$ and~$G$ can be defined over~$K$.
\end{definition}

\begin{remark}
We observe that independence is an equivalence
relation, since if $F(P_1)=G(P_2)$, then $(\hat G\circ F)(P_1)=(\hat
G\circ G)(P_2)$, where~$\hat G$ is the dual isogeny to~$G$.  We also
note that a torsion point can never be part of an independent pair,
since if (say) $NP_1=0$, then we can take $F=[N]$ and $G=0$ to show that~$P_1$
and any~$P_2$ are dependent.
\end{remark}

\begin{conjecture}
\label{conjecture:ECoverFF}
Let~$K$ be a characteristic zero function field as above, let~$E_1/K$
and~$E_2/K$ be elliptic curves, and let $P_1\in E(K)$ and
$P_2\in E_2(K)$ be $K$-independent points. (See
Definition~\ref{definition:independence}.) 
\begin{parts}
\Part{(a)}
There is a constant
$c=c(K,E_1,E_2,P_1,P_2)$ so that
\[
  \deg\GCD(n_1P_1,n_2P_2) \le c
  \qquad\text{for all $n_1,n_2\ge1$.}
\]
\Part{(b)}
Further, there is an equality
\[
  \GCD(nP_1,nP_2) = \GCD(P_1,P_2) 
  \qquad\text{for infinitely many $n\ge1$.}
\]
\end{parts}
\end{conjecture}

We prove Conjecture~\ref{conjecture:ECoverFF} in the case that~$E_1$
and~$E_2$ have constant $j$-invariant. We note that even this
``special case'' is far from trivial, since it relies on Raynaud's
Theorem (Theorem~\ref{theorem:raynaud}). For~(b), we prove a stronger
positive density result.

\begin{theorem}
\label{theorem:ECoverFFconstantj}
Let~$K$ be a characteristic zero function field as above, let~$E_1/K$
and~$E_2/K$ be elliptic curves, and let $P_1\in E(K)$ and $P_2\in
E_2(K)$ be $K$-independent points.  Assume further that the elliptic
curves $E_1/K$ and~$E_2/K$ both have constant $j$-invariant, i.e.,
\text{$j(E_1),j(E_2)\in k$}.
\begin{parts}
\Part{(a)}
There is a constant
$c=c(K,E_1,E_2,P_1,P_2)$ so that
\[
  \deg\GCD(n_1P_1,n_2P_2) \le c
  \qquad\text{for all $n_1,n_2\ge1$.}
\]
\Part{(b)}
The set
\[
  \bigl\{n\ge1 : \GCD(nP_1,nP_2) = \GCD(P_1,P_2)  \bigr\}
\]
has positive density.
\end{parts}
\end{theorem}
\begin{proof}
The fact that~$E_1$ and~$E_2$ have constant $j$-invariants means that
they split over some finite extension of~$K$. Taking a common
splitting field, there is a finite cover $C'\to C$ and elliptic
curves~$E_1',E_2'/k$ so that
\[
  \Ecal_i\times_CC'  \cong_{/k}
  E_i'\times_kC'.
\]
(N.B., $E_i'$ is defined over the constant field~$k$.)
We thus get commutative diagrams 
\[
  \begin{CD}
    E_i' \times_k C' @>>> \Ecal_i \\
    @VVV @VVV \\
    C' @>f>> C \\
  \end{CD}
  \qquad\text{for $i=1,2$.}
\]
Each point $P_i\in E_i(K)$ gives a section $\s_{P_i}:C\to\Ecal_i$,
which in turn lifts to a unique section 
\[
  \t_{P_i}\times 1 :C'\to E_i'\times_kC'. 
\]
In other words, each point $P_i\in E_i(K)$ gives a unique morphism
\text{$\t_{P_i}:C'\to E_i'$} so that the following diagram commutes:
\begin{equation}
\label{cd:sectionsplits}
\begin{CD}
  E_i'\times_k C' @>>> \Ecal_i \\
  @A\t_{P_i}\times 1AA @AA\s_{P_i}A \\
  C' @>f>> C \\
\end{CD}
\end{equation}
\par
We now fix two $K$-independent points $P_1\in E_1(K)$ and
$P_2\in E_2(K)$ and define a morphism
\[
  \f = \t_{P_1}\times\t_{P_2} : C' \longrightarrow E_1'\times_k E_2'.
\]
Suppose that a point $\g\in C$ is in the support of $\GCD(n_1P_1,n_2P_2)$
for some $n_1,n_2\ge1$. This means that
\[
  \g\in\Support(\s_{n_1P_1}^*(\Obar_{\Ecal_1}))
  \qquad\text{and}\qquad
  \g\in\Support(\s_{n_2P_2}^*(\Obar_{\Ecal_2})).
\]
Tracing around the commutative diagrams, this means that
for every point $\g'\in f^{-1}(\g)\in C'$,
\[
  \t_{n_1P_1}(\g')=O_1
  \qquad\text{and}\qquad
  \t_{n_2P_2}(\g')=O_2,
\]
where~$O_i\in E_i'(k)$ is the zero point.
Equivalently, $\t_{P_1}(\g')\in E_1'[n_1]$ and $\t_{P_2}(\g')\in E_2'[n_2]$,
so in particular,~$\t_{P_1}(\g')$ and~$\t_{P_2}(\g')$ are
torsion points of~$E_1'$ and~$E_2'$, respectively.
Hence $\f(\g')=\bigl(\t_{P_1}(\g'),\t_{P_2}(\g')\bigr)$ is a torsion point
of the abelian surface~\text{$E_1'\times E_2'$}.
\par
To recapitulate, we have proven that
\begin{multline}
  \label{equation:supportGCDistorsion}
  \g\in\Support(\GCD(n_1P_1,n_2P_2)) \\
  \Longrightarrow
  \f(\g')\in E_1'[n_1]\times E_2'[n_2]
  \subset (E_1'\times E_2')_\tors \\
  \quad\text{for all $\g'\in f^{-1}(\g)$.}
\end{multline}
To ease notation, we let
\[
  A=E_1'\times E_2' \qquad\text{and}\qquad V = \f(C')\subset A.
\]
There are several cases to consider:

\subsection*{Case I: $\t_{P_1}$ and $\t_{P_2}$ are both constant maps}
\hfill\break From diagram~\eqref{cd:sectionsplits}, the assumption
that~$\t_{P_i}$ is constant and nonzero implies that the divisor
$\s_{P_i}^*(\Obar_{\Ecal_i})$ is supported on the set of ramification
points~$\Rcal_f$ of the map \text{$f:C'\to C$}. 
More generally, the independence
assumption implies that~$P_i$ is nontorsion, so $nP_i\ne O_i$ for all
$n\ge1$. Hence $\s_{nP_i}^*(\Obar)$ is supported on~$\Rcal_f$.
Further, Lemma~\ref{lemma:stabilityofvanishing}(b) says that for any
particular point~$\g\in\Rcal_f$, the multiplicity
$\ord_\g\s_{nP_i}^*(\Obar_{\Ecal_i})$ is bounded independently
of~$n$. Therefore
\[
  \text{$\deg\s_{nP_i}^*(\Obar_{\Ecal_i})$ is bounded for all $n\ge1$.}
\]
Thus~$E_i(K)$ contains infinitely many points of bounded degree, or
what amounts to the same thing, of bounded height.
(See~\cite[III~\S4]{ATAEC}, where $h(P)=2\deg\s_P^*(\Obar)+O(1)$.)  It
follows from~\cite[Theorem~III.5.4]{ATAEC} that $\Ecal_i\to C$ splits
as a product over~$k$.\footnote{For Case~I, it is also possible to
give a more elementary proof that~$\Ecal_i$ splits using explicit
Weierstrass equations and considering different types of twists.}
\par
Thus Case~I leads to the conclusion that both~$E_1$ and~$E_2$
are $K$-isomorphic to elliptic curves defined over~$k$, so we may
replace them with curves that are defined over~$k$.  Then
$\Ecal_i=E_i\times_kC$, and any point~$Q_i\in E_i(K)$ is associated to a 
$k$-morphism $\t_{Q_i}:C\to E_i$. Our assumption that~$\t_{P_i}$ is 
constant is equivalent to saying that $P_i\in E_i(k)$, so as long
as $nP_i\ne O_i$, we have
\[
  \Support\left(\s_{nP_i}^*(\Obar_{\Ecal_i})\right) = 
  \bigl(\{nP_i\}\times C'\bigr) \cap   \bigl(\{O_i\}\times C'\bigr) 
  = \emptyset.
\]
Hence the assumption that~$P_1$ and~$P_2$ are nontorsion points leads,
in Case~I, to the conclusion that
\[
  \GCD(n_1P_1,n_2P_2) = 0\qquad\text{for all $n_1,n_2\ge1$.}
\]
Thus Case~I gives a strong form of both~(a) and~(b).

\subsection*{Case II: $\t_{P_1}$ or $\t_{P_2}$ is nonconstant,
and $V\cap A_\tors$ is infinite}
\hfill\break
The assumption that one of~$\t_{P_1}$ or~$\t_{P_2}$ is nonconstant implies
that~$V=\f(C')$ is an irreducible curve, and then Raynaud's
Theorem~\ref{theorem:raynaud} tells us that $V\cap A_\tors$ can only
be infinite if it is contained in the translate of an elliptic curve
(abelian subvariety of~$A$) by a torsion point of~$A$.  Thus there is
an elliptic curve $W\subset A$ and a torsion point $t\in A$ so that
$V=W+t$. Let~$N$ be the order of the point~$t$. Then composing with
the multiplication-by-$N$ map yields
\[
  [N]\circ\f = [N]\circ(\t_{P_1}\times\t_{P_2}) = \t_{NP_1}\times\t_{NP_2},  
\]
and since~$W$ is an elliptic curve, we see that~$[N]\circ\f$ maps~$C'$
onto \text{$NV=N(W+t)=NW=W$}. Hence we get a commutative diagram
\[
  \begin{matrix}
  & & C' \\
  & \overset{\t_{NP_1}}{\swarrow} 
  & \phantom{XX}\big\downarrow{\scriptstyle [N]\circ\f} 
  & \overset{\t_{NP_2}}{\searrow} \\[1ex]
  E_1' & \overset{\pi_1}{\longleftarrow} & W
      & \overset{\pi_2}{\longrightarrow} & E_2' \\
  \end{matrix}
\]
where~$\pi_1$ and~$\pi_2$ are the projections 
\text{$\pi_i:E_1'\times E_2'\to E_i'$}.
\par
Let $d_2=\deg(\pi_2)$.  Since~$W$ is an elliptic curve, there is a
dual isogeny $\hat\pi_2:E_2'\to W$ with the property that
$\hat\pi_2\circ\pi_2=[d_2]$. We compute
\begin{align}
  [d_2]\circ \t_{NP_1}
  &= [d_2]\circ \pi_1 \circ [N]\circ\f  \notag \\ 
  &= \pi_1 \circ [d_2] \circ [N]\circ\f  \notag \\
  &= \pi_1 \circ \hat\pi_2 \circ \pi_2 \circ [N]\circ\f  \notag \\ 
  &= \pi_1 \circ \hat\pi_2 \circ \t_{NP_2}
    \label{equation:P2P1relation}
\end{align}
\par
Let $G':E_2'\to E_1'$ be the isogeny
\[
  G' = \pi_1\circ\hat\pi_2\circ[N] \in \Hom_k(E_2',E_1').
\]
Recall that $K'=k(C')$ is the extension of~$K$ over which~$E_1$ and~$E_2$ 
become isomorphic to~$E_1'$ and~$E_2'$, respectively. Thus~$G'$ induces
an isogeny
\[
  G:E_2\to E_1\quad\text{defined over $K'$,}
\]
but \textit{a priori}, there is no reason that~$G$ need be defined over~$K$.
However, the relation~\eqref{equation:P2P1relation} gives a commutative
diagram
\[
  \begin{CD}
  E_2' \times_k C' @>G'\times1>> E_1'\times_k C' \\
  @A\t_{P_2}\times1AA  @AA\t_{[d_1N]P_1}\times1A \\
  C' @= C' \\
  \end{CD}
\]
which is equivalent to the equality
\begin{equation}
  \label{equation:Gmapsrationalpts}
  G(P_2) = [d_1N](P_1)
\end{equation}
of points in~$E_1(K)$.
\par
The curves~$E_1$ and~$E_2$ and the points~$P_1$ and~$P_2$ are rational
over~$K$ by assumption, hence the same is true of the multiple
$[d_1N](P_1)$ of~$P_1$.  Thus~\eqref{equation:Gmapsrationalpts} says
that the isogeny~$G$ maps at least one $K$-rational point of~$E_2$ to a
$K$-rational point of~$E_1$. Further, the independence assumption
on~$P_1$ and~$P_2$ ensures that they are not torsion points. Hence
Lemma~\ref{lemma:isogenyofrationalpt} tells us that~$G$ is indeed
defined over~$K$.  Then~\eqref{equation:Gmapsrationalpts} contradicts
the $K$-independence of~$P_1$ and~$P_2$, which shows that Case~II
cannot occur.  

\subsection*{Case III: $\t_{P_1}$ or $\t_{P_2}$  is nonconstant,
and $V\cap A_\tors$ is finite}
\hfill\break
The assumption that one of~$\t_{P_1}$ or~$\t_{P_2}$ is nonconstant implies
that the map~$\f$ is nonconstant, and hence that $\f:C'\to V$ is
finite-to-one. We showed earlier~\eqref{equation:supportGCDistorsion} that
\[
  \Support(\GCD(n_1P_1,n_2P_2)) \subset f\bigl(\f^{-1}(V\cap A_\tors)\bigr),
\]
so the assumption that $V\cap A_\tors$ is finite implies that
$\GCD(n_1P_1,n_2P_2)$ is supported on a finite set of points that is
\textit{independent of~$n_1$ and~$n_2$}. Since
Lemma~\ref{lemma:stabilityofvanishing}(b) tells us that for any
particular point~$\g\in C$, the order of~$\GCD(n_1P_1,n_2P_2)$ at~$\g$ is
bounded independently of~$n_1$ and~$n_2$, this shows that
$\deg\GCD(n_1P_1,n_2P_2)$ is bounded, which completes the proof of~(a).
\par
In order to prove~(b), we return to~\eqref{equation:supportGCDistorsion},
which actually provides the more accurate information that
\[
  \Support(\GCD(nP_1,nP_2)) \subset f\bigl(\f^{-1}(V\cap A[n])\bigr).
\]
Since $V\cap A_\tors$ is finite by assumption, we can find an integer~$N$
so that $V\cap A_\tors$ is contained in~$A[N]$. It follows that
\[
  V\cap A[n] = V\cap A[\gcd(n,N)]
  \qquad\text{for all $n\ge1$,}
\]
and hence in particular that
\[
  V\cap A[n] = V\cap \{0\}
  \qquad\text{for all $n$ with $\gcd(n,N)=1$.}
\]
Hence
\begin{multline*}
  \Support(\GCD(nP_1,nP_2)) = \Support(\GCD(P_1,P_2))\\
  \qquad\text{for all $n$ with $\gcd(n,N)=1$.}
\end{multline*}
On the other hand, Lemma~\ref{lemma:stabilityofvanishing}(a) tells us
that the multiplicities of $\GCD(nP_1,nP_2)$ and $\GCD(P_1,P_2)$ are
the same at every point in the support of the latter. Therefore
\[
  \GCD(nP_1,nP_2)=\GCD(P_1,P_2)
  \qquad\text{for all $n$ with $\gcd(n,N)=1$,}
\]
which completes the proof of~(b), and with it the proof of
Theorem~\ref{theorem:ECoverFFconstantj}.
\end{proof}

\begin{remark}
One can easily formulate other variants of the common divisor problem
on algebraic groups. For example, let $a(T)\in\CC[T]$ be a nonconstant
polynomial, let $E/\CC(T)$ be an elliptic curve, and let $P\in E(\CC(T))$
be a nontorsion point. Then it is plausible to guess that there is
a constant $c=c(a,E,P)$ so that
\[
  \deg\gcd(D_{nP},a^m-1)\le c
  \qquad\text{for all $n_1,n_2\ge1$.}
\]
If~$E$ has constant $j$-invariant, one can probably prove that this is
true using a generalization of Raynaud's theorem to semiabelian
varieties~\cite{C-L,DP}. The situation over~$\QQ$ is somewhat more
complicated due to the different growth rates of~$D_{nP}$ and~$a^m$,
but Vojta's conjecture applied to the blowup of \text{$E\times\GG_m$}
at $(0,1)$ implies that for every $\e>0$ there exists a proper Zariski
closed subset $Z=Z(a,E,P,\e)$ of \text{$E\times\GG_m$} so that
\begin{multline*}
  \log\gcd(D_{nP},a^m-1)\le \e\max\{n^2,m\} \\
     \text{provided that $(nP,a^m)\notin Z$.}
\end{multline*}
See~\cite{SilvermanGCDinFGgps} for details.
\end{remark}

%%%%%%%%%%%%%%%%%%%%%%%%%%%%%%%%%%%%%%%%%%%%%%%%%%%%%%%%%%%%%%%%%%%%%%
\section{Common Divisors on Elliptic Curves over 
Characteristic~$p$ Function Fields}

We continue with the notation set in Section~\ref{section:preliminaries},
except that rather than working over a field of characteristic~0, 
we work instead over a finite field~$k=\FF_q$.
\par
Let $a(T),b(T)\in k[T]$ be multiplicatively independent polynomials.
As noted in the introduction, Ailon and Rudnick~\cite{AR} prove that
\text{$\gcd(a(T)^n-1,b(T)^n-1)$} is bounded for $n\ge1$ when~$k$ is a
field of characteristic zero, but the author~\cite{SilvermanGCDoverFF}
has shown that there is no analogous bound when~$k$ has
characteristic~$p$, even if the exponent~$n$ is subject to some
reasonable restrictions such as~$p\notdivide n$.
\par
It is natural to ask for a result similar to~\cite{SilvermanGCDoverFF}
for elliptic curves over~$\FF_q(T)$, as given in the following
conjecture.

\begin{conjecture}
\label{conjecture:gcdonECoverFF}
Let~$\FF_q$ be a finite field of characteristic~$p$, let~$E/\FF_q(T)$
be an elliptic curve, and let~$P,Q\in E(\FF_q(T))$ be nontorsion
points.  Then there is a constant $c=c(q,E,P,Q)>0$ so that
\begin{equation}
  \label{equation:gcdonECoverFF}
  \deg\GCD(nP,nQ) \ge cn
  \qquad\text{for infinitely many $n\ge1$ with $p\notdivide n$.}
\end{equation}
\end{conjecture}

\begin{remark}
It is tempting to conjecture a lower bound of the form~$cn^2$, since
the only obvious upper bound comes from $\deg D_{nP}\gg\ll n^2$, but
there is really no evidence either for or against the stronger bound.
\end{remark}

\begin{remark}
It is easy to prove~\eqref{equation:gcdonECoverFF} if one
allows~$p$ to divide~$n$.  To see this, factor $[p]=\f\circ\fhat:E\to
E$, where $\f:E^{(p)}\to E$ is the Frobenius map and~$\fhat$ its dual.
Let~$\Obar$ denote, as usual, the zero divisor on a model of~$E$
over~$\PP^1$, and let~$\Obar'$ similarly denote the zero divisor on a
model of~$E^{(p)}$.  Then $\f^*(\Obar)=p\Obar'$ and
$\fhat^*(\Obar')=\Obar+D$ for some effective divisor~$D$, which allows
us to estimate
\begin{align*}
  \deg\GCD(p^iP,p^iQ) 
  &=\deg\GCD\left(\s_P^*\circ\left.{\fhat^i}\right.^*\circ{\f^i}^*(\Obar),
            \s_Q^*\circ\left.{\fhat^i}\right.^*\circ{\f^i}^*(\Obar)\right) \\
  &= p^i\deg\GCD\left(\s_P^*\circ\left.{\fhat^i}\right.^*(\Obar'),
                \s_Q^*\circ\left.{\fhat^i}\right.^*(\Obar')\right), \\
  &\ge p^i\deg\GCD\left(\s_P^*(\Obar),\s_Q^*(\Obar)\right).
\end{align*}
Hence
\[
  \deg\GCD(nP,nQ) \ge n\cdot\deg\GCD(P,Q)
  \quad\text{for all $n=p^i$, $i=1,2,3,\ldots$.}
\]
\end{remark}

We prove a strong form of Conjecture~\ref{conjecture:gcdonECoverFF}
for elliptic curves with constant $j$-invariant.

\begin{theorem}
\label{theorem:gcdonECoverFFconstantj}
Let~$\FF_q$ be a finite field of characteristic~$p\ge5$, 
let $E/\FF_q(T)$ be an elliptic curve, let~$P,Q\in E(\FF_q(T))$ be
nontorsion points, and suppose that $j(E)\in\FF_q$.  Then 
\begin{multline*}
  \deg\GCD(nP,nQ) \ge n + O(\sqrt{n})\\
  \text{for infinitely many $n\ge1$ with $p\notdivide n$,} 
\end{multline*}
where the big-$O$ constant depends only on~$E/\FF_q(T)$.
\end{theorem}
\begin{proof}
For the moment, we take~$E/\FF_q(T)$ to be any elliptic curve, not
necessarily with constant~$j$ invariant, and we fix a (minimal)
Weierstrass equation for~$E$.  For each integer $N\ge1$, let
\[
  S_{q,N} = \{\pi\in\FF_q[T] : \text{$\pi$ is monic, irreducible, and
  $\deg\pi=N$}\}.
\]
Given any $\pi\in S_{q,N}$, we reduce~$E$ modulo~$\pi$ to obtain
an elliptic curve $\Etilde_\pi$ defined over the finite field
$\FF_\pi=\FF_q[T]/(\pi)$. The residue fields~$\FF_\pi\cong\FF_{q^N}$
associated to the various~$\pi$ are all isomorphic, but the elliptic
curves~$\Etilde_\pi$ for different primes need not (and generally will
not) be isomorphic. The Hasse estimate~\cite[V.1.1]{AEC} says that
\begin{equation}
  \label{equation:hasseestimate}
  n_\pi(E) =\#\Etilde_\pi(\FF_\pi)=q^N + 1 - a_\pi(E)
  \quad\text{with $|a_\pi(E)|\le2q^{N/2}$.}
\end{equation}
\par
Suppose now that $j(E)\in\FF_q$. For simplicity, we assume that
$j(E)\ne0,1728$. The other two cases, which can be handled similarly,
will be left for the reader.  This means that there is an elliptic
curve $E'/\FF_q$ so that~$E$ is a quadratic twist of~$E'$.
More prosaically, if~$E'$ is given by a Weierstrass equation
$y^2=x^3+ax+b$ with $a,b,\in\FF_q^*$, then~$E$ has a Weierstrass
equation~(cf.~\cite[X~\S5]{AEC})
\[
  E : y^2=x^3+\d^2ax+\d^3b
  \qquad\text{for some squarefree $\d\in\FF_q[T]$.}
\]
Replacing~$a,b$ by~$r^2a,r^3b$ and~$\d(T)$ by $r^{-1}\d(T)$ for an
appropriate~$r\in\FF_q^*$, we may assume that~$\d(T)$ is monic.
For now, we assume that~$\d(T)\ne1$, so~$E$ is a nontrivial
twist of~$E'$.
\par
For any~$\pi\in S_{q,N}$ with $\pi\notdivide\d$, the
curve~$\Etilde_\pi/\FF_\pi$ is isomorphic over~$\FF_\pi$ to
either~$E'/\FF_{q^N}$ or to the unique quadratic twist
of~$E'/\FF_{q^N}$. More precisely,~$\Etilde_\pi/\FF_\pi$ is
isomorphic over~$\FF_\pi$ to~$E'/\FF_{q^N}$ if~$\d$ is a square
in~$\FF_\pi$ and it is isomorphic to the twist if~$\d$ is not
a square in~$\FF_\pi$. It follows from this and from the standard proof
of~\eqref{equation:hasseestimate} using the action of Galois on
the Tate module~\cite[V.1.1]{AEC} that
\[
  a_\pi(E) = \LS{\d}{\pi}a_N(E'),
\]
where $a_N(E')=q^N+1-\#E'(\FF_{q^N})$ and where $\LS{\d}{\pi}$ is the
Legendre symbol. We divide the set of primes~$S_{q,N}$ into two
subsets,
\begin{align*}
  S_{q,N}^+(\d) &= \left\{ \pi\in S_{q,N} : \LS{\d}{\pi}=+1 \right\}, \\
  S_{q,N}^-(\d) &= \left\{ \pi\in S_{q,N} : \LS{\d}{\pi}=-1 \right\}.
\end{align*}
Then
\[
  n_\pi(E) = \begin{cases}
             q^N + 1 - a_N(E')&\text{for all $\pi\in S_{q,N}^+$,} \\
             q^N + 1 + a_N(E')&\text{for all $\pi\in S_{q,N}^-$.} \\
             \end{cases}
\]
\par
For a fixed~$\d$, the (quadratic) reciprocity law
for~$\FF_q[T]$~\cite[Theorem~3.5]{RosenFF} says that 
\[
  \LS{\d}{\pi} = (-1)^{\frac{q-1}{2}\cdot N\cdot\deg(\d)}\LS{\pi}{\d}.
\]
Notice that the power of~$-1$ depends only on the degree of~$\pi$.  It
follows that half the possible congruence classes for~$\pi$
modulo~$\d$ yield~$\LS{\d}{\pi}=+1$ and the other half
yield~$\LS{\d}{\pi}=-1$. Now Dirichlet's theorem
\cite[Theorem~4.8]{RosenFF} implies that
\[
  \#S_{q,N}^+ = \frac{q^N}{2N} + O\left(\frac{q^{N/2}}{N}\right)
  \qquad\text{and}\qquad
  \#S_{q,N}^- = \frac{q^N}{2N} + O\left(\frac{q^{N/2}}{N}\right).
\]
\par
Let $n=q^N+1-a_N(E')$. Then~$n=n_\pi(E)$ for every $\pi\in S_{q,N}^+$,
so~$n$ annihilates~$\Etilde_\pi(\FF_\pi)$, and hence $D_{nP}$ is
divisible by all of these primes. Since the same is true of~$nQ$, we
obtain the lower bound
\begin{multline*}
  \deg\GCD(nP,nQ) 
  \ge \sum_{\pi\in S_{q,N}^+} \deg(\pi) \\
  =\#S_{q,N}^+\cdot N 
  = \frac{1}{2}q^N + O(q^{N/2})
  = n + O(n^{1/2}).
\end{multline*}
Similarly, if $n=q^N+1-a_N(E')$, then the same argument using
the primes $\pi\in S_{q,N}^-$ yields the same lower bound.
\par
To recapitulate, we have proven that
\begin{multline}
  \deg\GCD(nP,nQ)\ge n+O(n^{1/2}) \\
  \text{for all $n=q^N+1\pm a_N(E')$ with $N=1,2,3,\ldots$,}
  \label{equation:FFlowerbound}
\end{multline}
where we may take either sign. This estimate is exactly the lower
bound that we are trying to prove, subject to the additional constraint
that we want~$n$ to be relatively prime to~$p$. However, it is clear
that at least one of the numbers \text{$q^N+1+a_N(E')$}
and  \text{$q^N+1-a_N(E')$} is prime to~$p$, since otherwise~$p$ would
divide their sum, and hence~$p=2$, contrary to assumption. 
Therefore~\eqref{equation:FFlowerbound} holds for infinitely many values
of~$n$ with~$p\notdivide n$, which completes the proof of 
Theorem~\ref{theorem:gcdonECoverFFconstantj}.
\par
It remains to consider that case that~$E$ is a trivial twist of~$E'$,
i.e., the case that~$E$ is~$\FF_q(T)$-isomorphic to a cruve defined
over~$\FF_q$. But then $E(\FF_q(T))=E'(\FF_q(T))=E'(\FF_q)$, since a
nonconstant point in~$E'(\FF_q(T))$ would correspond to a nonconstant
morphism $\PP^1\to E'$. But the group~$E'(\FF_q)$ is finite,
so~$E(\FF_q(T))$ has no nontorsion points and the statement of the
theorem is vacuously true.
\end{proof}

\begin{remark}
We continue with the notation from the proof of
Theorem~\ref{theorem:gcdonECoverFFconstantj}. It is well known that
$a_N(E')=\a^N+\b^N$, where~$\a$ and~$\b$ are the complex roots
of $X^2-a_1(E')X+q$. Thus for any particular~$E'$, one may find
more precise information about the values of~$n$ being used in the
statement of the theorem.
\end{remark}

%%%%%%%%%%%%%%%%%%%%%%%%%%%%%%%%%%%%%%%%%%%%%%%%%%%%%%%%%%%%%%%%%%%%%%

\end{document}